\documentclass{amsart}
\usepackage{amsmath,amssymb,amscd,amsthm}
\usepackage[pdftex]{graphicx}

\newcommand{\rsp}{\raisebox{0em}[2.4ex][1.5ex]{\rule{0em}{2ex} }}

\title{Leonardo da Vinci's Proof of the Theorem of Pythagoras}
\author{Franz Lemmermeyer}
\begin{document}
\maketitle

\pagestyle{myheadings}
\markboth{Leonardo da Vinci's Proof}{\today \hfil Franz Lemmermeyer}

While collecting various proofs of the Pythagorean Theorem for
presenting them in my class (see \cite{Lemm}) I discovered a beautiful 
proof credited to Leonardo da Vinci. It is based on the diagram
on the right, and I leave the pleasure of reconstructing the
simple proof from this diagram to the reader (see, however, the
proof given at the end of this article).

\medskip\noindent
\begin{minipage}{7cm}
Since I had decided to give correct references to as many results 
as possible (if only to set an example) I started looking at the many 
presentations of da Vinci's proofs on the internet for finding 
out where da Vinci had published his proof. It turned out that 
although this was a very well known proof, none of the many 
sources was able to point me to a book with a sound reference
let alone to one of da Vinci's original publications. For this 
reason I gave up and simple remarked
\end{minipage} \quad
\begin{minipage}{5cm}
\includegraphics[width=5cm]{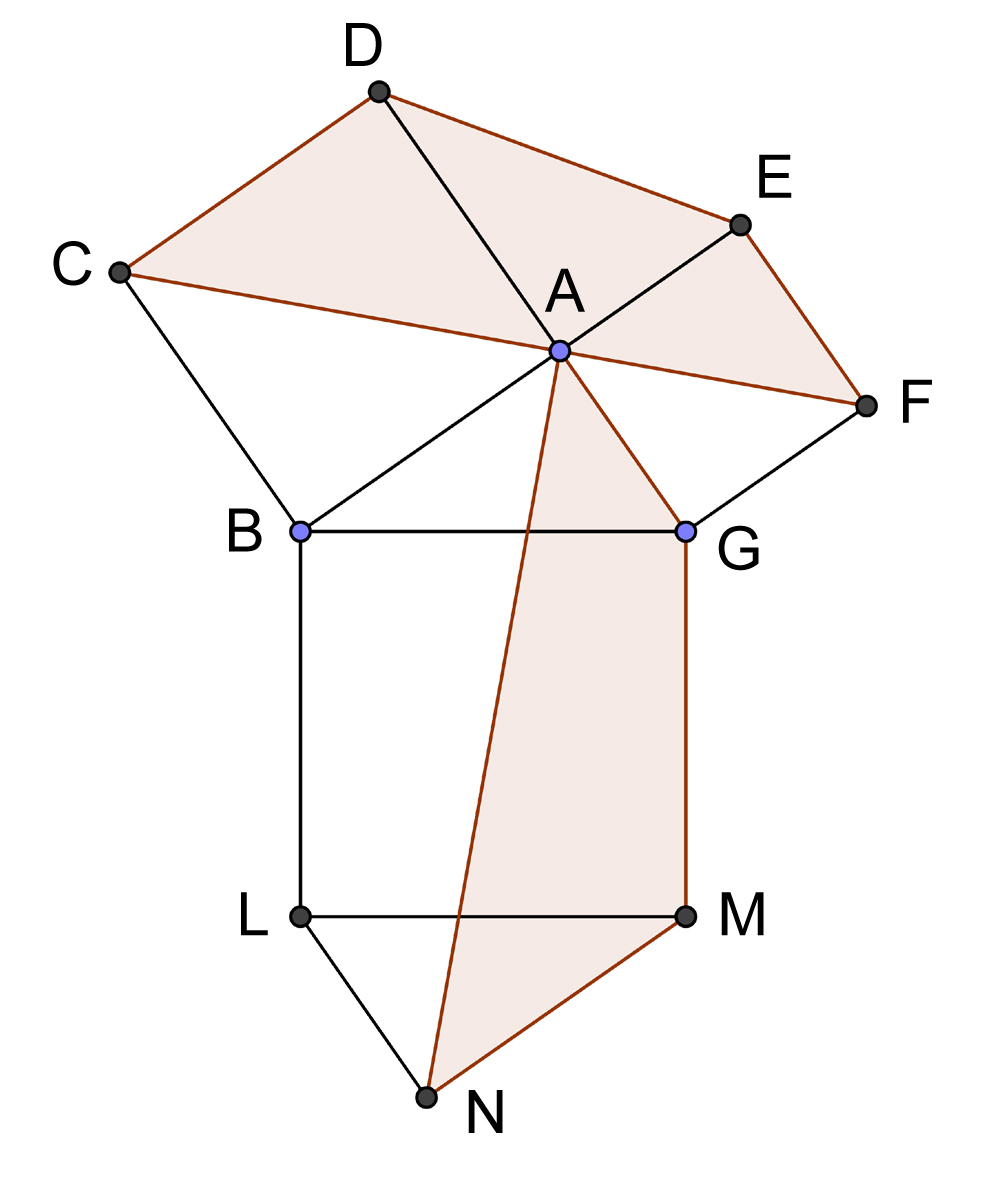}
\end{minipage}

\medskip\begin{quote}
{\em As a rule, one has to be very careful in such situations,}
\end{quote}
meaning that claims that cannot be verified often turn out to 
be wrong. A simple search (with the keywords ``da Vinci'' and 
Pythagoras or Pitagoras) will reveal dozens of books that credit
the proof above to da Vinci, including Maor's book \cite{MaorP} or 
the otherwise very accurate and highly readable book \cite{OW} by 
Ostermann \& Wanner, where the proof is credited to da Vinci in 
Exercise 21 on p.~26.

\subsection*{Tracking down References}

Eventually I developed an interest in finding the precise reference,
and I began hunting for the origin of da Vinci's proof. The first thing
I disovered, more or less by accident while browsing through old
geometry texts on the internet, was that a proof identical to da Vinci's 
proof was found by Terquem in 1838. In \cite[S.~103--104]{Terquem}  
he writes
\begin{quote}
{\em The theorem of Pythagoras being very important, we will give
     here a new proof based only on the superposition of figures.}
\end{quote}

Proofs of the Pythagorean Theorem have been rediscovered over and
over again, so the fact that Terquem had found a proof credited to 
da Vinci does not mean that da Vinci did not find it first. Terquem's
proof was republished in 1893 by M.~Balitrand \cite{JME}, without
any reference to da Vinci.

Next I found a remark in Eli Maor's wonderful book on the Pythagorean
theorem; in \cite[p.104]{MaorP} he writes
\begin{quote}
{\em Loomis, on the authority of F.C.~Boon, A.C. (Miscellaneous Mathematics,
1924) attributes this proof to Leonardo da Vinci (1452--1519).}
\end{quote}
Here the reference is to the first edition of Loomis \cite{Loomis}.
I was unable to find the publication by Boon. Boon did  publish a book 
in 1924, namely `` A companion to elementary school mathematics'', which 
is also referred to by Loomis. As Brian Hopkins informed me, the pages 
of this book actually do carry the title ``Miscellaneous Mathematics'', 
and it does credit the proof to da Vinci, suggesting that these two 
books actually are one and the same. Then I discovered \cite{DP}, 
whose authors also give the reference to Boon but mention in addition 
that the credit to da Vinci may be found in Heath's book \cite{HeathE}
(Boon also refers to Heath). Heath, who was extremely well versed in 
classical mathematics, writes in \cite[S.~365]{HeathE} 
\begin{quote}
{\em It appears to come from one of the scientific papers of Leonardo 
    da Vinci (1452--1519).}
\end{quote}
In addition he remarks that this proof may be found in compilations 
of proofs of the Pythagorean Theorem by J.W.~M\"uller and Ign.~Hoffmann.
Another search for these names led me to the report on elementary geometry
written by Max~Simon \cite{Simon}, where the reference to da Vinci is
also given. Indeed Heath quotes Simon very often, mainly his edition
of the first six books of the Elements, but also (on pp. 202, 328) 
his report \cite{Simon}.

Using these references I was finally able to track down what I think is
the correct historical development behind ``da Vinci's proof''.

\subsection*{First Appearance of ``da Vinci's'' Proof}

The story begins in 1790; in his book \cite[p.~124--126]{Tempelhof}
on ``Geometry for soldiers, and those who are not'', 
Tempelhof (some sources, including wikipedia, spell his name Tempelhoff; 
in his book, Tempelhof is used) gives the proof in question and remarks
\begin{quote}
{\em This proof is somewhat roundabout; but it has the advantage that
     the truth of the theorem can be made clearly visible.}
\end{quote}
Tempelhof does not give any information about who discovered the proof.

In 1819, Hoffmann \cite{HoffP} published a compilation of 32 proofs of the 
theorem of Pythagoras, using as his main sources earlier compilations, in 
particular the dissertations of Scherz \& St\"ober \cite{Stoe} 
(Strasbourg 1743; Scherz was St\"ober's supervisor -- his name figures
prominently on the cover of this dissertation) and by Lange \& Jetze 
\cite{Jetze} (Halle 1752; again, Lange was Jetze's supervisor). In the 
same year, M\"uller \cite[p.~64]{MuellSZ}
published, as a response to Hoffmann's publications, his own
compilation, which contained ``da Vinci's'' proof as no.~15 (out of 18
different proofs). In this regard he wrote (\cite[p.~64]{MuellSZ})
\begin{quote}
{\em I have known this construction only for a few years now since it 
     has been communicated to me orally.}
\end{quote}

This proof was one out of three that Hoffmann added to his 32 proofs in
the second edition \cite{HoffP2} of \cite{HoffP} along with the following 
comment:
\begin{quote}
{\em In this appendix three proofs are given. The first seems to have 
     been known for quite a while since it can be found in several 
     older writings. Its discoverer is not named.}
\end{quote}
I do not know to which books the ``older writings'' refer; all of the 
ancient books and dissertations in the bibliography are available through 
google books. It is a pity that the list of figures, which are usually 
presented on the last couple of pages in these geometry books, have not 
been reproduced properly. If these figures could be accessed, searching
for geometry textbooks containing the proof in question would be a lot
easier.

The discoverer was finally revealed in \cite[p. 70--71]{MuellNB} (1826) 
by M\"uller, a professor of mathematics at the gymnasium in Nuremberg:

\begin{quote}
{\em In this regard I remark that Mr. Hofrath Joh. Tobias Mayer from 
  G\"ottingen is the discoverer of the proof that I have
  given as the fifteenth in my book mentioned above, ``systematic 
  compilation'' etc. page 62 until 64. He has found this proof already
  in 1772 and has repeatedly presented it in his lectures in Altdorf
  given in the years 1779--1785 and has disseminated in this way. Therefore
  Tempelhof could include it in his Geometry for Soldiers published in 1790.}
\end{quote}
This sounds as if Mayer had seen the proof in M\"uller's book
and told him that this proof was due to himself.

I still do not know how Tempelhof learned about Mayer's proof. 
Tempelhof studied in Frankfurt (Oder) and Halle, and joined the 
Prussian army in 1756. He taught officers in Berlin as well as the
King's son and his brother, but apparently did not leave Berlin 
except for taking part in various military campaigns. Altdorf, on
the other hand, is a small town outside of Nuremberg; its university 
was closed in 1809.

\subsection*{Tobias Mayer}

For most of today's mathematicians, the name Tobias Mayer will be 
pretty much unknown. Johann Tobias Mayer was born in G\"ottingen in 1752
and died there in 1830. His father, Tobias Mayer (Marbach 1723 -- 
G\"ottingen 1762), published his first book on mathematics at the 
age of $18$. He died very young, and was famous for his accomplishments 
in astronomy -- he even received a part of the Longitude Prize
(see Forbes \cite{Forbes}) for his contribution to improving navigation 
on sea through his lunar tables. One of his inventions also was used
in the measurements of the arc of the meridian in connection with 
determining the metre. Mayer's correspondence with Euler was translated 
into English and published by Forbes.

Actually Tobias Mayer is mentioned in Gauss's letter to
Olbers from October 26, 1802:
\begin{quote}
{\em I do not know any professor who has done much for science except 
     the great Tobias Mayer, and he was regarded as a bad professor in 
     his times.}
\end{quote}
There is, by the way, also a connection between Gau\ss{} and
Tempelhof: Gau\ss{} had read both of Tempelhof's books on analysis, 
and even sent him a copy of his dissertation. In a letter to Bolyai 
from December 16, 1799, he writes that, in his opinion, General von 
Tempelhof is one of the best German mathematicians (in 1799, France
was the leading nation in mathematics, as is testified by mathematicians
such as Fourier, Laplace, Lagrange, Legendre, Monge, Poisson, Poinsot, 
and Poncelet. On the German side there was Gauss).

Mayer's son Johann Tobias Mayer began lecturing in G\"ottingen in 1773,
and moved to Altdorf in 1780; in 1786 he went to Erlangen, and in
1799 he returned to G\"ottingen. He is the author of various textbooks
on practical geometry and differential calculus, and published many
articles on physics. There is a very active Tobias-Mayer-Society and a 
Tobias-Mayer museum in Marbach.

\subsection*{Subsequent Development} 

The authorship of Johann Tobias Mayer is also mentioned subsequently by
Hoffmann in various of his writings; in his comments on Euclid's Elements 
\cite[p. 284]{HoffE} he writes
\begin{quote}
{\em Another highly simple and astute proof, whose discovery
is credited to Joh. Tob. Mayer in G\"ottingen, is the following.} 
\end{quote}

A similar remark can be found in Hoffmann \cite[p.~8]{HoffB}:
\begin{quote}
%
{\em The Pythagorean Theorem according to Johann Tobias Mayer and a 
     variation of this proof.

1) In my memoir ``{\em The Pythagorean Theorem, equipped with two and thirty  
   proofs that are partially known and partially new}'', Mainz 1821, 2nd. ed.
   p.~36, and in my memoir ``{\em The Elements of Euclides} etc. Mainz 1829,
   p.~284, I have mentioned the astute proof whose discovery is credited
   (by Joh. Wolfg. M\"uller ``{\em New contributions to the theory of 
   parallels}'' etc., Augsburg and Leipzig 1826, p.~70--71) to the 
   very sagacious geometer Joh. Tobias Mayer in G\"ottingen.}
\end{quote}

\subsection*{Enter da Vinci}

Thus the proof of the theorem of  Pythagoras based on the  diagram above 
is due to Mayer (1772) and was rediscovered by Terquem (1838). The 
question remains where da Vinci enters the picture. He does so in
Max Simon's report \cite{Simon} ``On the development of elementary
geometry in the 19th century'' published by the German Union of
Mathematicians in 1906, in the same series in which e.g. Hilbert's
report on algebraic numbers had appeared. Originally, Simon's report 
had been written for the encyclopedia of mathematics, but eventually 
Felix Klein rejected the submission. In fact, in his preface Simon writes
\begin{quote}
{\em Thus when Mr. Klein finally rejected the review in its present 
     form then this happened mainly because he did not have any 
     scientific assistants at his disposal who would give all references 
     with bibliographic precision at each place where some work is 
     mentioned. In fact, the condition of the slips of papers 
     necessitated an extremely time-consuming correction.}
\end{quote}
And on p.~111 we can finally read
\begin{quote}
{\em The proof of the hexagon, which was included in very many textbooks, 
     such as Mehler, is not due to  T\'ed\'enat (Manuel), but may be 
     found in the second edition of Hoffmann as no.~33 from older w
     ritings (Lionardo da Vinci).}
\end{quote}

The need for correcting the references is clearly visible here: the
textbook by Mehler is probably the one by M\"uller we have already cited, 
and my guess is that the ``manuel'' by ``T\'ed\'enat'' is actually 
Terquem's book \cite{Terquem}; Lionardo, of course, should read
Leonardo. There is a geometry textbook written by T\'edenat in the
year 7 of the revolution (1799), namely \cite{Ted}; there T\'edenat
gives two proofs of the theorem of Pythagoras, the first in art.~176
(the proof given by Euclid) and the second, a variation of Euclid's
proof, in art.~215.

At the end of Simon's memoir there are two pages listing 
all kinds of misprints, but his list of errata is far from being 
complete. Actually it is difficult to explain Klein's suggestion of 
publishing Simon's memoir ``in the present form'' as a report
in the same series of the Jahresberichte as Hilbert's report.

Simon takes the attribution for this proof from \cite{HoffP2}, and his 
expression ``older  writings'' ( ``\"altere Schriften'' in the German
original) is exactly the expression used by Hoffmann. Apparently Simon 
was quoting from memory and must have mixed up the ``older writings'' 
with another reference involving Leonardo da Vinci (this is the fourth 
out of four references to da Vinci in his report); another explanation 
might be that some ``scientific assistant'' trying to make sense of 
Simon's collection of slips of papers is responsible for this mix-up.

In any case my conclusion is that the legend of Leonardo da Vinci's proof 
was given birth by Max Simon's report from 1906, from where it was copied 
by Heath, Boon and Loomis and then was spread throughout the mathematical 
literature; nowadays the claim that Mayer's proof was found by da Vinci
can be found in dozens of books, hundreds of articles and thousands
of web pages.

\subsection*{Summary}

As a conclusion, here is a short summary of these developments:

\bigskip
\begin{tabular}{|r|p{10cm}|} \hline 
\rsp 1790 & Tempelhof gives ``da Vinci's'' proof in 
            \cite[p.~124--126]{Tempelhof} without attribution. \\
\rsp 1819 & M\"uller \cite[p.~64]{MuellSZ} gives the proof  
            without attribution. \\
\rsp 1821 & Hoffmann \cite[S.~36]{HoffP2} gives the proof and writes that
            its discoverer is unknown. \\
\rsp 1826 & M\"uller \cite[p. 70--71]{MuellNB} credits the proof to
            Johannes Tobias Mayer from G\"ottingen. \\
\rsp 1838 & Terquem \cite{Terquem} rediscovers Mayer's proof. \\
\rsp 1848 & Hoffmann \cite[p.~8]{HoffB} also credits the proof to Mayer, 
            referring to M\"uller. \\
\rsp 1906 & Max~Simon \cite{Simon} credits the proof to da Vinci for the 
            first time and refers to Hoffmann \cite{HoffP2}. \\
\rsp 1908 & Heath \cite{HeathE} reads about ``da Vinci's proof''
            in Simon's report, but also refers to M\"uller and Hoffmann. \\
\rsp 1924 & Boon \cite{Boon} credits the proof to da Vinci; in the 
            bibliography, Heath \cite{HeathE} is cited. \\
\rsp 1940 & Loomis \cite{Loomis} credits the proof to da Vinci on the 
            authority of Boon. \\
\hline
\end{tabular}
\bigskip

Despite my findings I still would love to hear from an expert on
da Vinci's mathematical writings that a proof of the Pythagorean 
theorem by da Vinci's hand is not known.

\section*{A Variation of Mayer's Proof}

The following version of Mayer's proof is perhaps a little bit
simpler than the proofs given elsewhere. Take the diagram with
the right-angled triangle and the squares on its sides, and 
rotate the square on the hypotenus along with the triangle by
$180^\circ$ around the center $Z$ of the square. By construction,
the are of the quadrilateral AGMN is $\frac12c^2 + F$, where $F$
is the area of the right-angled triangle.

\smallskip
\begin{minipage}{6cm}
\includegraphics[width=5cm]{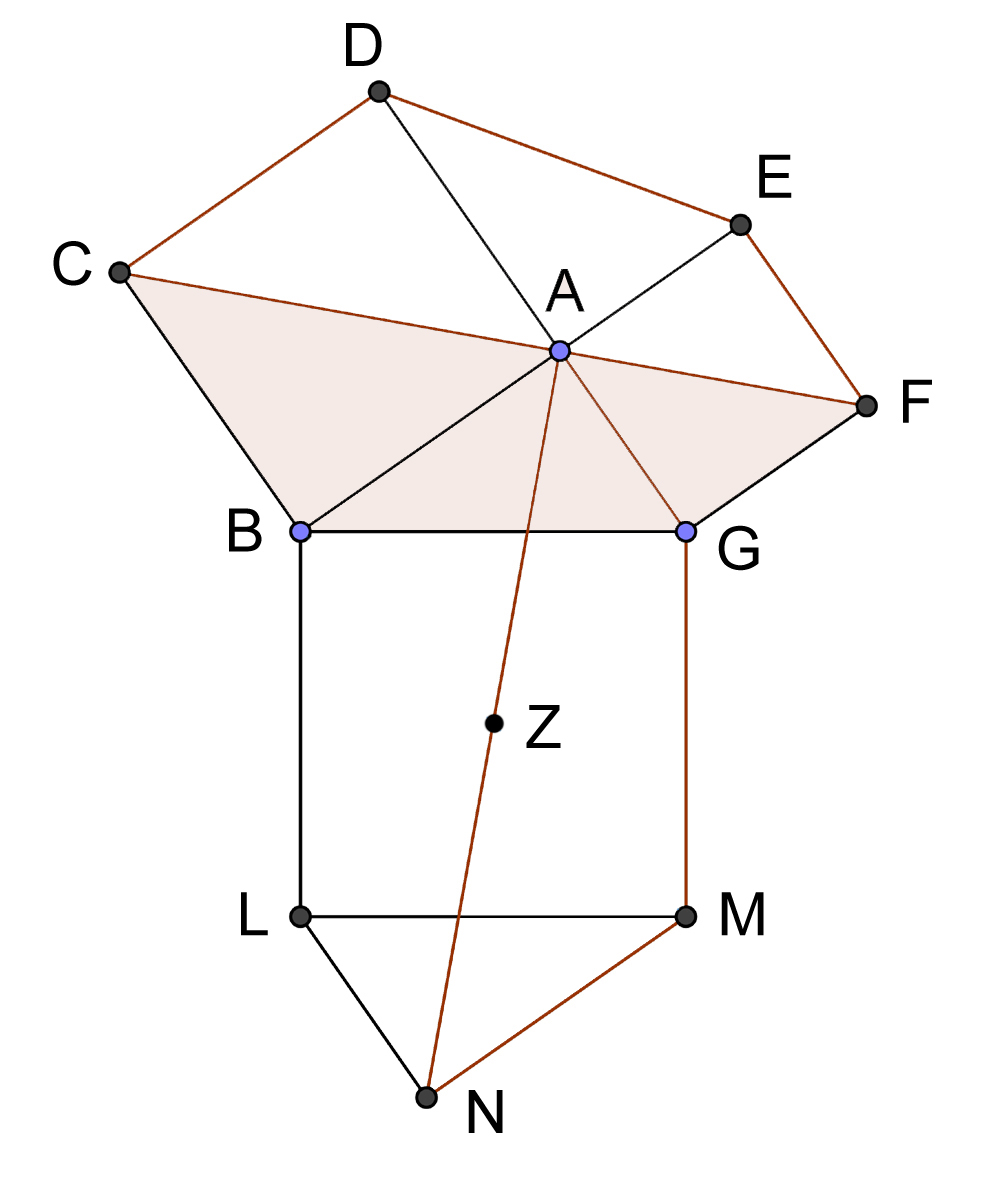}
\end{minipage} \quad \begin{minipage}{6cm}
\includegraphics[width=5cm]{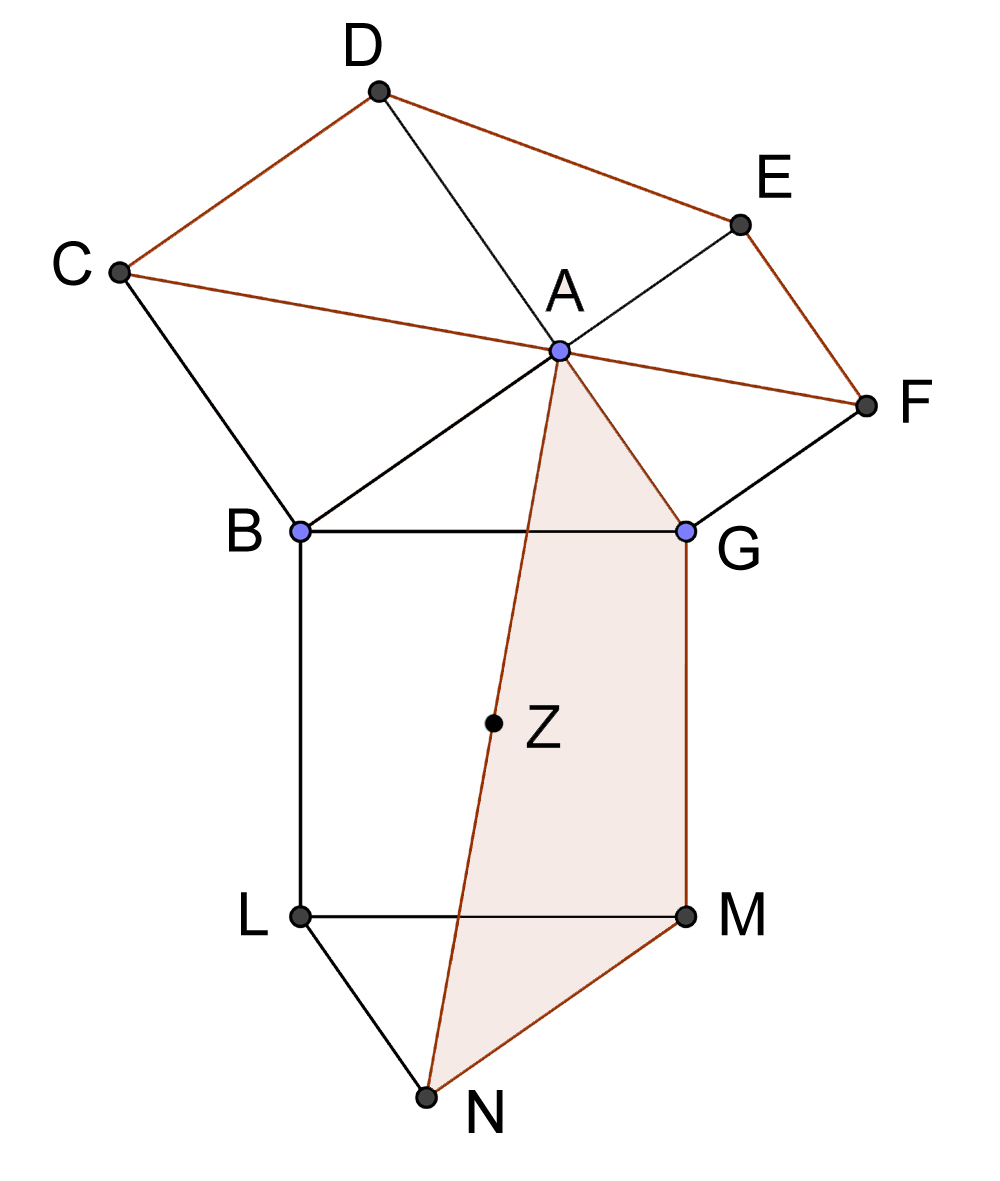}
\end{minipage}
\smallskip

On the other hand, the area of the quadrilateral BCFG is
equal to that of CFED since they are symmetric with respect 
to the line CF. Observe that A  is on this line since 
$\angle CAF = \angle CAB + \angle BAG + \angle GAF
            = 45^\circ + 90^\circ + 45^\circ = 180^\circ$.
Thus  BCFG has area $\frac12 a^2 + \frac12 b^2 + F$.

Finally, rotating BCFG by $90^\circ$ about G moves BCFG into AGMN,
so their areas are equal, and Pythagoras follows.

\section*{Acknowledgements}

I thank Brian~Hopkins for several suggestions and for providing me with 
many details concerning Boon's book.

\end{document}